\documentclass[12pt]{article}
\usepackage[margin=23mm]{geometry}                % See geometry.pdf to learn the layout options. There are lots.
\usepackage{amssymb,amsmath,amsthm}
\newtheorem{prop}{Proposition}
\usepackage{wrapfig}
\usepackage{graphicx,mathptmx}
\DeclareMathOperator{\K}{\mathit{C}\,}

\begin{document}

\title{Are Random Axioms Useful?}
\author{A.Shen}
\date{}
\maketitle
\let\eps=\varepsilon

The famous G\"odel incompleteness theorem says that for every sufficiently rich formal theory (containing formal arithmetic in some natural sense) there exist true unprovable statements. Such statements would be natural candidates for being added as axioms, but where can we obtain them? One classical (and well studied) approach is to add (to some theory $T$) an axiom that claims the consistency of $T$. In this note we discuss the other one (motivated by Chaitin's version of the G\"odel theorem) and show that it is not really useful (in the sense that it does not help us to prove new interesting theorems, see below Proposition~\ref{conservation}). However, the situation could change if we take into account the size of the proofs.

\section{G\"odel theorem in Chaitin's version}

Kolmogorov complexity $\K(x)$ of a binary string $x$ is defined as the minimal length of a program (without input) that outputs $x$ and terminates. (This definition depends on a programming language, and one should choose one that makes complexity minimal up to $O(1)$ additive term.) Most strings of length $n$ have complexity close to $n$. More precisely, the fraction of $n$-bit strings that have complexity less than $n-c$,  is at most $2^{-c}$. In particular, there exist strings of arbitrary high complexity. 

However, as G.~Chaitin pointed, the situation changes if we look for strings of \emph{provably} high complexity. More precisely, we are looking for strings $x$ and numbers $n$ such that the statement ``$\K(x)>n$'' (for these $x$ and $n$, so it is a closed statement) is provable in formal arithmetic. Chaitin noted that all $n$ that appear in these statements, are less than some constant $c$. (The argument is a version of the Berry paradox: assume for arbitrary integer $k$ we can find some string $x$ such that $\K(x)>k$ is provable; let $x_k$ be the first string in the order of enumeration of all proofs; this definition provides a program of size $O(\log k)$ that generates $x_k$, which is impossible for large $k$ if $\K(x_k)>k$.)

Another proof of the same result shows that Kolmogorov complexity is not very essential here; by a standard fixed-point argument one can construct a program $p$ (without input) such that for every program $q$ (without input) the assumption that $p$ is equivalent to $q$ (produces the same output if $p$ terminates, and does not terminate if $p$ does not terminate) is consistent with formal arithmetic. If $p$ has length $k$, for every $x$ we may assume without contradiction that $p$ produces $x$, so one cannot prove that $\K(x)>k$.  

\section{Incompressibility axioms}

This leads to a natural idea: take a random string $x$ of length $n$ and consider the statement $\K(x)\ge n-1000$. It is true unless we are extremely unlucky. The probability of being unlucky is less than $2^{-1000}$. In statistics and natural science we are accustomed to identify this with impossibility. So we can add this statement and be sure that it is true; if $n$ is large enough, we get a true non-provable statement and could use it as a new axiom. We can even repeat this procedure several times: if the number of iterations $m$ is not astronomically large, $2^{-1000}m$ is still astronomically small.

Now the question: \emph{can we obtain a richer theory in this way and get some interesting consequences, still being practically sure that they are true}? 

The answer that we give (using a very simple argument): 

\begin{itemize}
\item yes, this is a safe way of enriching our theory (formal arithmetic), see Proposition~\ref{correctness};
\item yes, we can get a stronger theory in this way (Chaitin's theorem), but
\item no, we cannot prove anything interesting in this way, see Proposition~\ref{conservation}.
\end{itemize}

\section{Random axioms: soundness}

Let us first prove that the chances to get a false statement in this way are small. In fact, we can allow even more flexible procedure of adding random axioms that does not mention Kolmogorov complexity explicitly. This procedure goes as follows: assume that we already have proved for some number (numeral) $N$, for some rational $\delta>0$ and for some property $R(x)$ (with one free string variable $x$) that \emph{the number of strings $x$ of length $N$ such that $\lnot R(x)$ does not exceed $\delta 2^N$} (this is a closed arithmetic formula). Then we are allowed to toss a fair coin $N$ times, generating a string $r$ of length $N$, and consider the formula $R(r)$ as a new axiom. But we have to pay $\delta$ for each operation of this type, and initially we have only some fixed amount $\eps$. Intuitively, $\eps$ measures the maximal probability that we agree to consider as ``negligible''.

In this framework we shall speak not about \emph{proofs} (as sequences of formulas constructed according to some rules like \emph{modus ponens} etc.), but about \emph{proof strategies} that say which steps should be made for every outcome (i.e., for every string $r$ if proof strategy comes to randomized step). For a proof strategy $\pi$, formula $F$ and initial capital $\eps$ we can measure the probability that $\pi$ achieves (proves) $F$. The following proposition says that this procedure indeed can be trusted:

\begin{prop}[soundness]\label{correctness}
If the probability to achieve $F$ for a proof strategy $\pi$ with initial capital $\eps$ is greater than $\eps$, then $F$ is true. 
\end{prop}

This proposition says that the probability of proving a false statement is bounded by $\eps$. Moreover, the probability of the event ``a false statement appears in the proof'' is bounded by~$\eps$. This is intuitively obvious because a false statement means that one of the bad events has happened, and the sum of probabilities is bounded by $\eps$. However, this argument cannot be understood literally since different branches have different bad events (after branching). Formally one should prove by backward induction that at every point of the tree where the remaining capital is $\gamma$ and there are still no false statements, the conditional probability to get a false statement assuming that we have reached this point, is bounded by $\gamma$. (End of proof.) 

\section{Conservative extension?}

As Chaitin theorem shows, there are proof strategies that with high probability lead to \emph{some} statements that are non-provable (in formal arithmetic). However, the situation changes if we want to get some \emph{fixed} statement, as the following proposition (which is a stronger version of Proposition~\ref{correctness}) shows:

\begin{prop}[conservation]\label{conservation}
If the probability to get $F$ for a proof strategy $\pi$ with initial capital $\eps$ is greater than $\eps$, then $F$ is provable.
\end{prop}

Let us start with a simple case where there is only one randomized step. So some formula $R(x)$ is fixed, as well as some integer $N$, and it is provable (in formal arithmetic without extensions) that \emph{there are at most $\eps 2^N$ strings $x$ of length $N$ such that $\lnot R(x)$}. We denote this statement by $(*)$ in the sequel. Let us say that a string $x$ of length $N$ is \emph{strong} if the formula $R(x)\Rightarrow F$ is provable in arithmetic. The assumption guarantees that the fraction of strong $x$ among strings of length $N$ is greater than $\eps$. Let $x_1,\ldots,x_t$ be all of them. Then $R(x_i)\to F$ is provable for each $i$, therefore
    $$ R(x_1)\lor R(x_2)\lor\ldots\lor R(x_t) \rightarrow F$$
is also provable. But the disjunction in the left hand side is provably weaker than $(*)$, and the latter is provable according to the assumption. Therefore, $F$ is also provable. 

In the general case we need to be a bit more careful. We say that a vertex $u$ in a proof strategy tree is \emph{strong} if $F$ provably follows from the statements accepted up to that point (on the way to $u$) in the proof process. For every  vertex $u$ in the proof tree we also consider the remaining capital $\rho(u)$, and the probability $p(u)$ to finish the proof successfully starting from that point (=the conditional probability of successful proof under the condition that it goes through $u$). Then we prove the following statement using backward induction (from leaves to the root) over $u$: if $p(u)>\rho(u)$, then the vertex $u$ is strong. For the tree root this is the statement of the Proposition; for leaves the value of $p(u)$ is either $0$ or $1$, and if it is $1$, then $u$ is strong. 

The induction step uses the same argument as in the previous paragraph. Consider some vertex $u$ of the proof tree. Let $\delta(u)$ be the threshold used for branching at $u$. If the fraction of strong vertices among the sons of $u$ is greater than $\delta(u)$, then the vertex $u$ itself is strong (since the disjunction of statements $R(x)\to F$ over all $x$ that correspond to strong sons, is provable in $u$). So let us assume that this fraction does not exceed $\delta(u)$. Then at $u$ the random choice will lead to a strong vertex with probability at most $\delta(u)$. And if this is not the case, the chances of finishing the proof successfully (for each weak son) do not exceed the remaining capital $\rho(u)-\delta(u)$ by induction assumption. So the total probability of finishing the proof successfully starting at $u$ does not exceed $\rho(u)$. (End of proof.)

\section{Polynomial size proofs}

The situation changes drastically if we are interested in the length of proofs. The argument used in Proposition~\ref{conservation} gives an exponentially long proof compared with the original probabilistic proof (since we need to combine the proofs for all terms in the disjunction). (Here the complexity of a probabilistic proof is measured as the total length of formulas in the branch where this total length is maximal.) May be one can find another construction that shows that our ``probabilistic'' proofs can be transformed into a normal ones with only polynomial increase in size? Some reason why this should be hard is provided by the following Proposition~\ref{complexity}.

\begin{prop}\label{complexity}
Assume that such a polynomially bounded transformation is possible. Then complexity classes PSPACE and NP coincide.
\end{prop}

It is enough to consider standard PSPACE-complete language, the language TQBF of true quantified Boolean formulas. A standard interactive proof for this language (see, e.g., \cite{sipser}) uses Arthur--Merlin scheme where the verifier (``Arthur'') is able to perform polynomial-size computations and obtain random bits (that cannot be corrupted by Merlin). The correctness of this scheme is based on simple properties of finite fields. This kind of proof can be easily transformed into a successful probabilistic proof strategy in our sense. To explain this transformation, let us recall some details of the construction. Assume that a quantified Boolean formula $F$ starting with an universal quantifier is given. This formula is transformed into a statement that says that for some polynomial $P(x)$ two values $P(0)$ and $P(1)$ are equal to $1$. This polynomial is implicitly defined by a sequence of operations. To convince Arthur that it is indeed the case, Merlin shows Arthur this polynomial (listing explicitly its coefficients; there are polynomially many of them). In other terms, Merlin notes that the formula $$[\forall x (P(x)=\overline{P}(x))]\Rightarrow [P(0)=1 \land P(1)=1]$$ and therefore $$[\forall x(P(x)=\overline{P}(x))]\Rightarrow F$$ is true (and provable in arithmetic). Here $P(x)$ is the polynomial $P$ defined as the result of the sequence of operations, while $\overline{P}(x)$ is the explicitly given expression for $P$. Then Merlin notes that the implication $$[P(r)=\overline{P}(r)]\Rightarrow \forall x(P(x)=\overline{P}(x))$$ is true for most elements $r$ of the finite field and this fact is provable (using basic results about finite fields), so such an implication can be added as a new axiom. Continuing in this way, we get a probabilistic proof strategy that mimics Merlin's behavior in the IP-protocol for TQBF.

This strategy gives probabilistic proofs of polynomial length (assuming our proof system is not unreasonable inefficient). If they could be transformed into non-probabilistic proofs (and the underlying theory is consistent), these non-probabilistic proofs would be NP-witnesses for TQBF, and therefore PSPACE=NP. (End of proof.)

\section{Adding full information about complexities}

Let us be more generous and add to formal arithmetic full information about Kolmogorov complexity of all strings; for each string $x$ we add an axiom ``$\K(x)=k$'' where $k$ is the literal representing Kolmogorov complexity of $x$. (Of course, this procedure is non-constructive.) What theory do we get? It turns out that this theory can be easily described in logical terms:

\begin{prop}\label{full-information}
   This theory can be obtained from formal arithmetic by adding all true universal statements.
\end{prop}

First, note that the statements $\K(x)>n$ (for some string $x$ and number $n$) are universal statements (all programs of length at most $n$ do not produce $x$). The statements of the form $\K(x)\le n$ are existential and therefore are provable when true. So adding all true universal statements is enough to get full information about complexities. 

More interesting is the another direction. Now we assume that we have all the information about Kolmogorov complexities; we need to show that every true universal statement becomes provable. This
can be done by an argument which is a logic translation of the Turing completeness result for the set of pairs $(x,n)$ such that $\K(x)<n$. Here is it.

The Kolmogorov complexity function is enumerable from above: simulating in parallel all the programs, we get an upper bound for $\K$, and this upper bound decreases with time. Let $T_n$ be the moment when this upper bound reaches its final value (i.e., the true value of $\K$) for all strings of length at most $n$. 

Note that \emph{$\K(T)\ge n-O(\log n)$ for every $T\ge T_n$.} Indeed, knowing $n$ and $T$, one can find the true values of $\K$ for all strings of length $n$, and choose the first one that has complexity at least $n$ (it always exists). So $n$ and $T$ are enough to construct a string of complexity $n$, and the complexity of $n$ is only $\log n$.

Another simple observation: every program of length $n-O(\log n)$ either terminates after $T_n$ steps or does not terminate at all. Indeed, each terminating program may be considered as a description of the number of steps it takes to terminate, so the complexity of this number does not exceed the length of the program, and it remains to use the previous statement.

Both observation can be formalized as provable statements of formal arithmetic. If we have (as additional axioms) full information about Kolmogorov complexities of all strings, then for each $n$ one can prove a statement $T_n=t$ for some numeral $t$. Having such a statement, we can then prove non-termination of every non-terminating program of size at most $n-O(\log n)$. And this is enough to prove every true universal statement, since it can be reformulated as non-termination of some program. (End of proof.)

\textbf{Remark}. A sceptic could (rightfully) claim that the exact value of Kolmogorov complexity can encode some additional irrelevant information (in particular, about universal statements' truth values). For example, it is not obvious \emph{a priori} that the theory considered in Proposition~\ref{full-information} does not depend on the choice of optimal programming language fixed in the definition of Kolmogorov complexity. Also one may ask whether we consider plain or prefix complexity. 

To address all these questions, let us consider a smaller set of axioms and assume that for every $x$ the complexity of $x$ is guaranteed up to a factor of $2$, i.e., some axiom $\K(x)>c_x$ is added where $c_x$ is at least half of the true complexity of $x$. (We do not assume anything else about the choice of $c_x$; this choice cannot be done computably for obvious reasons.) Still the argument works with minor modifications. Let $\K^t(x)$ be the upper bound for complexity obtained if we restrict computation time by $t$. Let $x(n,t)$ be the string of length at most $n$ that makes $\K^t(x)$ maximal (among all $x$ of length at most $n$, for given $t$). Then $\K^t(x(n,t))$ is at least $n$ though $\K(x(n,t))$ can be much smaller (for small~$t$). By $T_n$ let us denote the minimal $T$ such that $\K^{T_n}(x)\le 2\K(x)$ for all strings $x$ of length at most $n$. The following statement is true for every $t$ and $n$:
    $$
t>T_n \Rightarrow \K(x(n,t))\ge n/2. 
    $$
Indeed, $\K^{t}(x(n,t))\ge n$ by definition, and if $t>T_n$, this guarantees $\K(x(n,t))\ge n/2$. Since $x(\cdot,\cdot)$ is computable, $\K(x(n,t))\le \K(t)+O(\log n)$. Therefore,
   $$
t>T_n\Rightarrow \K(t)\ge n/2-O(\log n).
   $$
Again, this implies that every program of size at most $n/2-O(\log n)$ that does not terminate in $T_n$ steps, does not terminate at all. This is not only true, but also provable (with natural formal representation of all involved notions). And, having our additional axioms, we are able (for each $n$) find some provable (using the axioms) upper bound for $T_n$, and therefore prove non-termination for every non-terminating program of size at most $n/2-O(\log n)$.

\textbf{Question}: What happens if we add only some axioms about complexities? Is it possible, for example, to add axioms $\K(x)\ge n-O(1)$ for most of the strings $x$ of length $n$ in such a way that a statement $\varphi$ (fixed in advance) remains unprovable?

\section{Acknowledgements}

The notions and questions considered in this note are rather natural, so probably they were considered long ago (though I have no specific indications). They were discussed during a colloquium talk given by the author at St.~Petersburg Euler institute few years ago~\cite{euler2006}. Recently essentially the same question was asked by Richard Lipton in his blog~\cite{lipton2011} (thanks to I.~Razenshteyn for showing me this post). L.~Beklemishev suggested to write down the answers (despite of their simplicity).

\end{document}